\documentclass[11pt]{article}

\usepackage[english]{babel}
\usepackage[T1]{fontenc}
\usepackage[utf8]{inputenc}
\usepackage{amsmath,amsfonts,amsthm,graphicx,color,mdwlist}
\usepackage{qsymbols,booktabs,psfrag,color}
\usepackage[%
            linktocpage,%
            bookmarks=true,%
            colorlinks=true,%
            ]{hyperref}
\usepackage{textcomp,xspace}

\def \eref#1{(\ref{#1})}
\def \1{{\bf 1}}
\def \l{\left}
\def \r{\right}

\def \ben{\begin{eqnarray}}
\def \een{\end{eqnarray}}
\def \be{\begin{eqnarray*}}
\def \ee{\end{eqnarray*}}

\def \b{{\sf b}}

\def \beq{\begin{equation}}
\def \eq{\end{equation}}

\def \Mult{{\sf Multinomial}}
\def \Pois{{\sf Poisson}}

\def \Var{{\sf Var}}
\def \Cov{{\sf Cov}}

\newtheorem{lem}{Lemma}[section]
\newtheorem{theo}{Theorem}%[section]
\theoremstyle{definition}
\def \dd{\xrightarrow[n]{(d)}}

\newcommand{\Rea}{\mathbb{R}}

\def \sous#1#2{\mathrel{\mathop{\kern 0pt#1}\limits_{#2}}}
\def \sur#1#2{\mathrel{\mathop{\kern 0pt#1}\limits^{#2}}}

\DeclareMathOperator{\cov}{cov}
\newcommand{\Hyp}{\textsf{Hyp}}

\newqsymbol{`E}{\mathbb{E}}
 \newqsymbol{`P}{\mathbb{P}}
\newqsymbol{`R}{\mathbb{R}}
\newqsymbol{`e}{\varepsilon}
\newqsymbol{`o}{\omega}
\def \bfun{T}

\renewcommand{\baselinestretch}{1.2}
\begin{document}
\begin{center}

\LARGE{\bf A functional central limit theorem for the \\
partial sums of sorted i.i.d. random variables}\medskip\normalsize
\[\begin{array}{ll}
\textrm{\Large Jean-Fran\c{c}ois Marckert}& \textrm{\Large ~~~~~~~~~~David Renault}\end{array}\]
\textrm{CNRS, LaBRI, Universit\'e de Bordeaux}\\
\textrm{351 cours de la Lib\'eration}\\
\textrm{33405 Talence cedex, France}\\
\textrm{email: name@labri.fr}
 \end{center}

 \begin{abstract} Let $(X_i,i\geq 1)$ be a sequence of i.i.d. random variables with values in $[0,1]$, and $f$ be a function such that $`E(f(X_1)^2)<+\infty$. We show a functional central limit theorem for the process $t\mapsto \sum_{i=1}^n f(X_i)1_{X_i\leq t}$. 
 \end{abstract}

\section{Introduction}

Let $(X_1, X_2, \ldots)$ be a sequence of i.i.d. random variables (r.v.) with values in $[0,1]$, having distribution~$\mu$, distribution function $F$, and defined on a common probability space $(\Omega,{\cal A},`P)$ on which the expectation operator is denoted $`E$. In this paper we are interested in proving a functional limit theorem for the sequence of processes $(Z_n,n\geq 1)$ defined by 
\ben\label{eq:ZN}
Z_n(t):= \frac1n{\sum_{i=1}^n f(X_i) \1_{X_i\leq t}},~~~ t\in[0,1]
\een
where $f:[0,1]\to `R$ is a measurable function. Let
$(\widehat{X}_i,1\leq i\leq n)$ be the sequence $(X_i,1\leq i\leq n)$
sorted in increasing order, and for any $t\in[0,1]$, denote by 
\[N_n(t)=\#\{i~: 1\leq i\leq n, X_i\leq t\}\]
the number of $X_i$'s smaller than $t$. Clearly, for $t\in[0,1]$, 
\[Z_n(t)=\frac1n{\sum_{k=1}^{N_n(t)} f(\widehat{X}_i)}.\]
Hence, $Z_n$ encodes the partial sums of functions of sorted
i.i.d. r.v., as mentioned in the title of this paper. In
order to state a central limit theorem for $Z_n$ the existence of
$\Var(f(X_1))<+\infty$ is clearly needed, but it is not sufficient to
control the fluctuations of $Z_n$ on all intervals. Standard considerations
about the binomial distribution implies that $N_n(t_2)-N_n(t_1)$ is quite
concentrated around $n(F(t_2)-F(t_1))$ (for $t_1<t_2$). 
Conditionally on $(N_n(t_1),N_n(t_2))=(n_1,n_2)$,
\beq
Z_n(t_2)-Z_n(t_1)\sur{=}{(d)}\frac{1}{n}\sum_{k=1}^{n_2-n_1} f(X_{(t_1,t_2]}(k))\eq
where $\sur{=}{(d)}$ means ``equals in distribution'', and where $(X_{(t_1,t_2]}(k), 1\leq k \leq n_2-n_1)$ 
is a family of i.i.d. r.v., whose common distribution is that of $X$ conditional on $X\in(t_1,t_2]$.
Hence, to get a functional central limit theorem for $Z_n$, the variances of these distributions need to be controlled. The following hypothesis $\Hyp$ is designed for that purpose:
\begin{quotation}
\noindent $\Hyp$: there exists an increasing function
$\bfun:[0,1]\to\Rea^+$ such that:
\[\left\{
\begin{array}{l}
x/\bfun(x) \textrm{ is bounded},\\
\bfun(x) \ln(x) \underset{x\to  0}{\longrightarrow} 0, \\
\forall I~\textrm{interval}~\subset [0;1], \quad
\Var\l(f(X)\,|\, X\in I\r) \leq \displaystyle\frac{\bfun(\mu(I))}{\mu(I)}
\end{array}\right.
\]
where $\Var(g(X)\,|\, X\in I)$ denotes the variance of $g(X)$ conditional on $X \in I$ (by convention, we set $`E(g(X)\,|\, X\in
I)=0$ when $`P(X\in I)=0$).
\end{quotation} 
When $f$ is bounded by $\gamma$ on $[0;1]$, the function  $\bfun(x) = \gamma^2 x$ satisfies $\Hyp$ (see also the discussion below Theorem \ref{theo:main}).  \par

Consider the mean of $Z_n$
\ben Z(t):= `E(Z_n(t))= `E\l(f(X) 1_{X\leq t}\r),\een  
(this can be shown to be a càdlàg process when $`E(|f(X)|)<+\infty$) and 
\ben
Y_n(t)=\sqrt{n}\l[Z_n(t)-Z(t)\r].
\een

The aim of this paper is to show the following result~:
\begin{theo}\label{theo:main} Let $(X_i,i\geq 0)$ be a sequence of i.i.d. r.v. taking their values in $[0,1]$ and $f:[0,1]\to `R$ a measurable function satisfying $\Hyp$, then
\be
Y_n & \dd &Y
\ee 
in $D[0,1]$, the space of càdlàg functions on [0,1] equipped with the Skorokhod topology, where  $(Y_t,t\in[0,1])$ is a centered Gaussian process with variance function
\ben\label{eq:var}
\Var(Y_s)&=&F(s)\Var(f(X)\,|\,X\leq s)+F(s)(1-F(s))`E(f(X)\,|\,X\leq s)^2
\een
and with covariance function, for $0\leq s<t\leq 1$
\ben\label{eq:cov}
\Cov(Y_s,Y_t-Y_s)&=& -F(s)(F(t)-F(s)) E(f(X)\,|\,X\leq s)`E(f(X)\,|\,s< X\leq t).
\een
\end{theo}
We discuss a bit the conditions in the theorem. Assume that the $X_i$'s are i.i.d. uniform on $[0,1]$, and that $f(x)=1/x^\alpha$ for some $\alpha>0$. The r.v. $f(X)=1/X^\alpha$ possesses a variance iff $\alpha<1/2$, and then it is in the domain of attraction of the normal distribution only in this case (Theorem~\ref{theo:main} needs this hypothesis for the convergence of $Y_n(1)$). The largest $\Var(f(X)|X\in(a,a+`e))$ is obtained for $a=0$, in which case we get
\[\Var\l(f(X) \,|\, X\in (0,`e]\r)= \frac{`e^{-2\alpha} \alpha^2}{(1-2\alpha)(1-\alpha)^2},\]
and one can check that $\alpha<1/2$ is also the condition for the existence of a function $\bfun$ satisfying $\Hyp$. $\Hyp$ appears to be a minimal assumption in that sense. \medskip

The first result concerning the convergence of empirical processes is due to Donsker's Theorem~\cite{D54}. It says that when $f$ is constant equal to $1$, then $Y_n$ converges in $D[0,1]$ to the standard Brownian bridge $\b$ up to a time change. A kind of miracle arises then, since the same analysis works for all distributions $\mu$ by a simple time change. This is not the case here. \par
Apart from strong convergence theorems à la Komlós-Major-Tusnády \cite{KMT}, modern results about the convergence of empirical processes --  see Shorack \& Wellner \cite{SW} and van der Vaart \& Wellner \cite{VW} -- much rely on the concept of Donsker classes, which we discuss below.

Denote by $`P_n=\frac1n\sum_{k=1}^n \delta_{X_i}$ the empirical measure associated with the sample $(X_i,1\leq i\leq n)$. As a measure, $`P_n$ operates on any set ${\cal F}$ of measurable functions $\phi:[0,1]\to `R$, 
\[`P_n\phi=\int_x \phi(x)d`P_n(x)=\sum_{k=1}^n \phi(X_i)/n.\]
 The empirical process is the signed measure $\mathbb{G}_n:=\sqrt{n}(`P_n-\mu)$. 
 By the standard central limit theorem, for a given function $\phi$ (such that $\mu \phi^2 <+\infty$), $\mathbb{G}_n\phi\dd {\cal N}(0,\mu(\phi-\mu \phi)^2)$, where ${\cal N}(m,\sigma^2)$ designates the normal distribution with mean $m$ and variance $\sigma^2$.

A \textit{P-Donsker class} is a set of measurable functions ${\cal F}$ such that $\l(\mathbb{G}_n \phi, \phi\in {\cal F}\r)$ converges in distribution to $\l(\mathbb{G} \phi,\phi\in {\cal F}\r)$,  in the $L_\infty$ topology (it is a central limit theorem for a process index by a set of functions). This means that~:
\begin{itemize} 
\item the convergence of the finite dimensional distributions holds~: (meaning that for any $k$, any $\phi_1,\dots,\phi_k \in {\cal F}$, $\l(\mathbb{G}_n \phi_1,\ldots, \mathbb{G}_n \phi_k\r)\dd N:=(N_1,\ldots,N_k)$
and $N$ is a centered Gaussian vector with covariance matrix $\Cov(N_i,N_j)=\mu\l[(\phi_i-\mu \phi_i)(\phi_j-\mu \phi_j)\r]$.
\item the sequence $\l(\mathbb{G}_n \phi, \phi\in {\cal F}\r)$ is tight in $L_\infty$.
\end{itemize}
The proof that a set forms a Donsker class is usually not that simple, and  numerous criteria can be found in the literature. In our case, the set of functions ${\cal F}$ is the following one~:
\[{\cal F}_f=\{ (x\mapsto \phi_t(x)=f(x)1_{x\leq t}), t\in[0,1]\}.\]  
We were unable to find such a criterion for this class, but notice that if such a result existed, it would imply Theorem \ref{theo:main} only for the topology $L_\infty$, a topology which is weaker than ours. Of course, Theorem \ref{theo:main} implies that ${\cal F}_f$ forms a Donsker class.\medskip

\noindent{\bf Note.} In fact classes ${\cal F}_f$ for non decreasing $f$, or for functions $f$ whose level sets are given by two intervals at most (such that $x\mapsto x^2$, $x\mapsto \cos(2\pi x)$, $x\mapsto \sin(2\pi x)$) are Donsker, since they are VC subgraph class (see Vapnik \& Chervonenkis \cite{VC}).\footnote{We thank Emmanuel Rio for this information}

If we consider the variables $X_i$'s in the formula \eref{eq:ZN}, as random times, then $Z_n(t)$ corresponds (up to the normalisation) to the sum of $f(X_i)$ for all events $X_i$ appearing before time $t$, where $f$ is some cost function. The process $Y_n$ appears to be the suitable tool to measure the fluctuations of $Z_n$. \par
We would like to mention \cite{MR}, a work at the origin of the present paper, written by the same authors. In  \cite{MR}, the convergence of rescaled trajectories made with sorted increments (in $\mathbb{C}$) to a deterministic convex is shown. For this purpose a weaker version of Theorem \ref{theo:main} is established.\medskip

We provide a proof of our theorem in an old fashioned style. We prove the convergence of the finite dimensional distributions, and then establish the tightness in $D[0,1]$; even if the proof is a bit technical, we think that several tricks make it interesting in its own right.

\section{Proof of Theorem \ref{theo:main}}

The proof starts with that of the convergence of the finite dimensional distributions (FDD) convergence of $Y_n$: this is classical as we will see. Let  $\theta_0:=0< \theta_1<\theta_2<\cdots<\theta_K=1$ for some $K\geq 1$ be fixed. In the sequel, for any function (random or not) $L$ indexed by $\theta$, $\Delta L(\theta_j)$ will stand for $L(\theta_j)-L(\theta_{j-1})$.  For any $\ell\leq
K$
\begin{equation}\label{eq:W}
\Delta Y_n(\theta_\ell)=\sqrt{n} \l[\Delta Z_n(N_n(\theta_j))- \Delta Z(\theta_j)\r],
\end{equation}
where by convention $Z_n(N_n(\theta_{-1}))=Z(\theta_{-1})=0$. The convergence of the FDD of $Y_n$ follows the convergence in distribution of the increments $(\Delta Y_n(\theta_\ell), 0\leq \ell \leq K)$. Notice that 
\beq
\Delta Z(\theta_j)= `E\l(f(X) 1_{\theta_{j-1}<X \leq \theta_j}\r).
\eq
If for some $j$, $\theta_{j-1}$ and $\theta_j$ are chosen in such a way that $\Delta F(\theta_j)=0$ then the $j$th increment in \eref{eq:W} is 0 almost surely (this is the case for the 0th increment if $\mu(\{0\})=0$). We now discuss the asymptotic behaviour of the other increments~: let $J=\{j \in \{0,\dots,K\} : \Delta F(\theta_j)\neq 0\}$.

Let $(n_j,j \in J)$ be some fixed integers summing to $n$. Denote by $\mu_{\theta_{j-1},\theta_j}$ the law of~$X$ conditioned by $\{\theta_{j-1}<X\leq \theta_j\}$. Conditional on $(N_n(\theta_j)=n_j,j \in J)$, the variables  $\Delta Z_n(N_n(\theta_j))$, $j\in J$  are independent, and  $\Delta Z_n(N_n(\theta_j))$ is a sum of $n_j-n_{j-1}$ i.i.d. copies of variables under $\mu_{\theta_{j-1},\theta_j}$, denoted from now on $(X_{\theta_{j-1},\theta_j}(k),k\geq 1)$. 

Since $\l(\Delta N_n(\theta_j),j \in J\r)\sim \Mult\l(n, (\Delta F(\theta_{j}), j \in J)\r)$, 
\beq\label{eq:limit_N}
 \l(\frac{\Delta N_n(\theta_j)- n\Delta F(\theta_{j})}{\sqrt{n}} ,j \in J\r)\dd (G_j,j\in J)\eq
where $(G_j,j \in J)$ is a centered Gaussian vector with covariance function, 
\[\cov(G_k,G_\ell)=-\Delta F(\theta_{k})\,.\, \Delta F(\theta_{\ell})+1_{k=\ell}\Delta F(\theta_k),\]
formula valid for any $0\leq k,\ell\leq K$. Putting together the previous considerations, we have
\ben\label{eq:deltaY_n}
\Delta Y_n(\theta_j)&=& \sum_{m=1}^{\Delta N_n(\theta_j)} \frac{f(X_{\theta_{j-1},\theta_j}(m))-`E(f(X_{\theta_{j-1},\theta_j}))}{\sqrt{n}}\\\label{eq:deltaY_n2}
                    & &+\l(  \frac{\Delta N_n(\theta_j)- n\Delta F(\theta_j)}{\sqrt{n}}\r) `E(f(X_{\theta_{j-1},\theta_j}))
\een
Using \eref{eq:limit_N} and the central limit theorem, we then get that
\beq\label{eq:FDD}
(\Delta Y_n(\theta_j),0\leq j \leq K)\dd  \l(\sqrt{\Delta F(\theta_j)}\widetilde{G}_j+G_j `E(f(X_{\theta_{j-1},\theta_j})),0\leq j \leq K\r)
 \eq
where the family of r.v. $(G_j,j\leq K)$ and $(\tilde{G_j},j\leq K)$ are independent, and the r.v. $\tilde{G_j}$ are independent centered Gaussian r.v. with variance $\Var(f(X_{\theta_{j-1},\theta_j}))$ (this allows one to determine the variance and covariance \eref{eq:var} and \eref{eq:cov}). Notice that here only the finiteness of $\Var(f(X_{\theta_{j-1},\theta_j}))$ and $`E(f(X_{\theta_{j-1},\theta_j}))$ are used. \medskip

It remains to show the tightness of the sequence $\l(Y_n,n\geq 0\r)$ in $D[0,1]$. A criterion for the tightness in $D[0,1]$ can be found in Billingsley \cite[Thm. 13.2]{BIL}: a sequence of processes $(Y_n,n\geq 1)$ with values in $D[0,1]$ is tight if, for any $`e\in(0,1)$, 
 \[\lim_{\delta\to 0} \limsup_n`P(`o'(Y_n,\delta)\geq `e)=0\]
where $`o'(f,\delta)=\inf_{(t_i)}\max_i \sup_{s,t\in[t_{i-1},t_i)} |f(s)-f(t)|,$ and the partitions $(t_i)$ range over all partitions of the form $0=t_0<t_1<\cdots<t_n\leq 1$ with $\min\{t_i-t_{i-1}, 1\leq i\leq n\}\geq \delta$. 

We now compare our current model formed by a set $\{X_1,\dots,X_n\}$ of $n$ i.i.d. copies of $X$ denoted from now on by $`P_n$, with a Poisson point process $P_n$ on $[0,1]$ with intensity $n \mu$, denoted by $`P_{P_n}$. Conditionally on $\#P_n=k$, the $k$ points
 $P_n:=\{X_1',\dots,X_k'\}$ are i.i.d. and have distribution $\mu$, and then $`P_{P_n}(~\cdot~ |\# P=n)=`P_n$.  The Poisson point process is naturally equipped with a filtration $\sigma:=\l\{\sigma_t=\sigma(\{ P\cap[0,t]\}),t\in[0,1]\r\}$.

We are here working under $`P_{P_n}$, and we let $N(\theta)=\#( P_n \cap [0,\theta])$; notice that under $`P_n$, $N$ and $N_n$ coincide.\par
Before starting, recall that if $N\sim \Pois(b)$, for any positive $\lambda$,
 \ben\label{eq:c-p1}
`P(N\geq x) &=&`P(e^{\lambda N}\geq e^{\lambda x})\leq `E(e^{\lambda N-\lambda x})
=e^{-b+be^{\lambda}-\lambda x}\\
\label{eq:c-p2}`P(N\leq x)&=&`P(e^{-\lambda N}\geq e^{-\lambda x}) \leq `E(e^{-\lambda N+\lambda x})=e^{-b+be^{-\lambda}+\lambda x}.
\een 

We explain now why the tightness of $(Y_n,n\geq 1)$ under $`P_{P_n}$ implies the same result under~$`P_n$. Let  $m=\inf\{x \in [0,1], F(x)\geq 1/2\}$ be  the median of $\mu$. 
\begin{lem} \label{lem:decoupe}There exists a constant $\gamma$ (which depends on $\mu$), such that for any $\sigma_m$-measurable event $A$,
\beq\label{eq:inter}
`P_n(A)=`P_{P_n}(A \,|\, \#P=n) \leq \gamma\, `P_{P_n}(A).
\eq 
\end{lem}
\noindent{\bf Proof of the Lemma} We have
\be
`P_{P_n}(A \,|\, \#P=n) &=&\sum_{k} \frac{`P_{P_n}(A ,\#(P\cap[0,m])=k)`P(\#P\cap[m,1]=n-k)}{`P(\#P=n)}\\
         & \leq &\sum_{k} `P_{P_n}(A ,\#(P\cap[0,m])=k)\sup_{k'}\frac{`P(\#P\cap[m,1]=n-k')}{`P(\#P=n)}\\
         & \leq & \gamma\,`P_{P_n}(A)
\ee
where $\gamma=\sup_{n\geq 1} \sup_{k'}\frac{`P(\#P\cap[m,1]=n-k')}{`P(\#P=n)}$, which is indeed finite since  $`P(\#P=n)\sim  (2\pi n) ^{-1/2}$, and since $\#P\cap[m,1]\sim \Pois(n/2)$, and then the probability that its value is $k$ is bounded above by some $d/\sqrt{n}$  according to Petrov \cite[Thm. 7 p. 48]{PET}. ~$\Box$ \medskip

Thanks to Lemma \ref{lem:decoupe}, if the sequence of restrictions $(Y_n|[0,m],n\geq 1)$ of $Y_n$ on $[0,m]$ is tight on $D[0,m]$ under $`P_{P_n}$ then so it is under $`P_n$ (the same proof works on $D[m,1]$ by a time reversal argument). 
To end the proof, we show that $(Y_n|[0,m],n\geq 1)$ is indeed tight under $`P_{P_n}$.

Take then some (small) $\eta\in(0,1)$, $`e>0$; we will show that one can find a finite partition $(t_{i},i \in I)$ of $[0,m]$ and a  $\delta\in(0,m)$ such that 
\beq\label{eq:tightness}
\limsup_n`P_n(`o'(Y_n,\delta)\geq `e)\leq \eta,
\eq
which is sufficient for our purpose. 

We decompose the process $Y_n$ as suggested by \eref{eq:deltaY_n} and \eref{eq:deltaY_n2},
\ben  
Y_n(\theta)&=& Y_n'(\theta)+Y_n''(\theta)
\een
where 
\ben\label{eq:Y_n'}
 Y_n'(\theta)&=&\sum_{m=1}^{N_n(\theta)} \frac{f(X_{[0,\theta]}(m))-`E(f(X_{[0,\theta]}))}{\sqrt{n}}\\
 \label{eq:Y_n''}Y_n''(\theta)&=& \l(  \frac{N_n(\theta)- n F(\theta)}{\sqrt{n} }\r) \frac{Z_\theta}{F(\theta)}.
\een
(If $F(\theta)$ then set $Y_n''(\theta)=0$ instead of \eref{eq:Y_n''}).\par
The tightness of each of the sequences $(Y_n',n\geq 1)$ and $(Y_n'',n\geq 1)$ in $D[0,1]$ suffices to deduce that of $(Y_n,n\geq 1)$. We then proceed separately. 

\subsection*{Tightness of $(Y_n',n\geq 1)$}
To control the jumps of $Y_n'$, we will need to localise the large atoms of $\mu$. 
Let $A=\{x \in [0,m], \mu(\{x\})>0\}$ be the set of positions of the atoms of $\mu$ in $[0,m]$, and let $A^{\ge a}:=\{x \in A : \mu(\{x\})\geq a\}$. Clearly $\#A^{\ge a}\leq 1/a$ and $[0,m]\setminus A^{\geq a}$ forms a finite union of open connected intervals $(O_{x},x\in G)$, with extremities $(t'_i,i \in I)$. The intervals $(O_x,x\in G)$ can be further cut as follows:\\
-- do nothing to those such that $\mu(O_x)<2a$,\\  
-- those such that $\mu(O_x)>2a$ are further split. Since they contain no atom with mass $>a$, they can be split into smaller intervals having all their weights in $[a,2a]$ except for at most one (in each interval $O_x$ which may have a weight smaller than $a$). 

Once all these splittings have been done, a list of at most $3/a$ intervals are obtained (in fact less than that), all of them having a weight smaller than $2a$. Name $G_a=(O_x, x\in I_a)$ the collection of obtained open intervals, indexed by some set $I_a$, and by $(t_i^a, i\geq 0)$ the partitions obtained. Take $O$ one of these intervals. One has $\#(P_n\cap O)$ is Poisson with parameter $n\mu(O)\leq na$.  Consider again \eref{eq:deltaY_n}, \eref{eq:deltaY_n2} and $\Hyp$. Set, for any $L\geq 1$,
\[S_L^{(n)}:=\sum_{\ell=1}^{L} \frac{f(X_{O}(\ell))-`E(f(X_{O}))}{\sqrt{n}}.\]
 Let 
\[`o(Y_n',O)=\sup\{|Y_n'(s)-Y_n'(t)|,s,t\in O\}\] 
be the modulus of continuity of $Y_n'$ on $O$. We have, for any $\alpha \in(0,1/2)$,  
\ben\label{eq:azz}
`P(`o(Y_n',O)\geq x) &\leq& `P\l( \l|{\#(P_n\cap O)- n
  \mu(O)}\r|\geq n^{1/2+\alpha} \r)\notag \\
                               &+&\sup_{L \in  \Gamma_n}
 `P\l( \sup\l\{\l|S_i^{(n)}-S_j^{(n)}\r|, i,j \leq L\r\}\geq x\r)
\een
where
\[\Gamma_n=\left[ n \mu(O)- n^{1/2+\alpha},
n\mu(O)+n^{1/2+\alpha} \right].\]
Using \eref{eq:c-p1} and \eref{eq:c-p2}, one sees that
\[`P\l( |P(n\mu(O))-n\mu(O)|  \geq n^{\alpha+1/2} \r) \leq c e^{-c'n^{\alpha}}\]
for some $c>0,c'>0$ and $n$ large enough (for this take $x=n\mu(O) +
n^{1/2+\alpha}$, $\lambda=1/\sqrt{n}$ in \eref{eq:c-p1} and,
$x=n\mu(O) -n^{1/2+\alpha}$, $\lambda=1/\sqrt{n}$ in \eref{eq:c-p2}).

Let us take care of the second term in \eref{eq:azz}. Clearly,
\be
 \sup\l\{\l|S_i^{(n)}-S_j^{(n)}\r|, i,j \leq L\r\}&=&\max_{i\leq L} S_i^{(n)}-\min_{j\leq L} S_j^{(n)}.
\ee
According to Petrov \cite[Thm.12 p50]{PET},
\be
`P\l(\max_{i\leq L} S_i^{(n)}\geq x\r) & \leq & 2`P\l(S_L^{(n)}\geq x -\sqrt{\frac{2 L \Var(f(X_{O}))}{n}}\r),
\ee 
and then 
\be
`P\l(\max_{i\leq L} S_i^{(n)}\geq x\r) & \leq & 2`P\l(S_L^{(n)}\geq x
-C_n(O)\r),
\ee
for $C_n(O)=\sqrt{\frac{2 L \bfun(\mu(O))}{n \mu(O)}}$,
and a similar inequality holds for  $\min_{i\leq L} S_i^{(n)}$. Since
\ben
`P\l(\max_{i\leq L} S_i^{(n)}-\min_{j\leq L} S_j^{(n)}\geq x\r)&\leq& `P\l(\max_{i\leq L} S_i^{(n)}\geq x/2\r)+ `P\l(-\min_{j\leq L} S_j^{(n)}\geq x/2\r)\notag\\
\label{eq:bound2}&\leq& 2`P\l(S_L^{(n)}\geq \frac{x}{2} -C_n(O)\r)+2`P\l(S_L^{(n)}\leq -\frac{x}{2} +C_n(O)\r).\notag
\een
To get some bounds, we use the central limit theorem for $S_L^{(n)}$,
and take $x=`e$, $a>0$ such that $\bfun(a)=`e^2\delta^2$ for some small
$\delta>0$ (recall that $\bfun$ is increasing and therefore
invertible), and any sequence $L_n$ such that $L_n/n\to \mu(O)$ (any sequence $L=L_n$ such that $L_n\in \Gamma_n$ satisfies this, and then we can control the supremum with this method).  We have 
\be
`P\l( S_L^{(n)}\geq \frac{`e}{2}- C_n(O)\r)&=&
`P\l(\frac{S_L^{(n)}}{\sqrt{ \mu(O) \Var(f(X_{O}))}}\geq \frac{`e/2 - C_n(O)}{\sqrt{\mu(O) \Var(f(X_{O}))}}\r).
\ee
For $n$ large enough,
\[C_n(O)\leq \sqrt{4 \bfun(\mu(O))}
\leq 2 `e\delta\]
and therefore 
\be
\limsup_n`P\l( S_L^{(n)}\geq \frac{`e}{2}- C_n(O)\r) 
&\leq & {\Phi}\l(\frac{`e/2 - 2 `e\delta}{\sqrt{\mu(O) \Var(f(X_{O}))}}\r)
\ee
where $\Phi$ is the tail function of the standard Gaussian distribution.

Finally, if $\delta$ is chosen sufficiently small ($2\delta <
1/2$), since $\mu(O)\Var(f(X_{O}))\leq \bfun(\mu(O)) \leq \bfun(a) =
`e^2\delta^2$, then on each interval $O\in G_a$,
\[`P\l(\sup\l\{\l|S_i^{(n)}-S_j^{(n)}\r|, i,j \leq L\r\} \geq `e\r)\leq 4\Phi\l(\frac{1/2 - 2\delta}{\delta}\r)\] 
and this independently of the choice of the interval $O$ in $G_a$, for $n$ large enough.   \par

The control of the intervals all together can be achieved using the union bound : since they are at most $3/\bfun^{-1}(`e^2\delta^2)$ such intervals, by the union bound 
\[`P_{P_n}\l(\sup_{O\in  G_a} `o( Y_n', O)\geq `e\r)\leq \frac{3}{\bfun^{-1}(`e^{2}\delta^2)}\l(4\Phi\l(\frac{1/2 - 2\delta}{\delta}\r)+c e^{-c'n^{\alpha}}\r).\]
Since $\Phi(x)\sous{\sim}{x\to+\infty} \exp(-x^2/2)/(\sqrt{2\pi}x)$,
and $\bfun(x) \ln(x) \underset{x\to 0}{\longrightarrow} 0$, which
implies that for any $`e>0$, and $\gamma > 0$ there exists a $\delta$
sufficiently small such that

\[\bfun(e^{-\gamma/\delta^2}) < `e^2 \delta^2 \quad \textrm{or equivalently} \quad
\frac{1}{\bfun^{-1} (`e^2\delta^2)} < e^{\gamma/\delta^2}\]
and as a result the probability can be taken as small as wanted. ~$\Box$

\subsection*{Tightness of $(Y_n'',n\geq 1)$}
Recall \eref{eq:Y_n''}. We work here under $`P_n$ and we only consider
the interval $I=\{\theta : F(\theta)>0\}$ since $Y_n''(\theta)$ equals
0 on its complement.  Since on $I$,
$\theta\mapsto\frac{Z_\theta}{F(\theta)}$ is càdlàg (and does not
depend on~$n$), it suffices to see why $\l(\frac{N_n(\theta)- n
  F(\theta)}{\sqrt{n}},n\geq 0\r)$ is tight in $D[0,1]$, but this is
clear since this is a consequence of the convergence of the standard
empirical process (Donsker \cite{D54}). ~$\Box$

\end{document}